\documentclass[12pt,reqno]{amsart}
\usepackage{amsmath}
\usepackage{amsxtra}
\usepackage{amscd}
\usepackage{amsthm}
\usepackage{amsfonts}
\usepackage{amssymb}
\usepackage{eucal}

\newtheorem{theorem}{Theorem}[section]
\newtheorem{conj}[theorem]{Conjecture}
\newtheorem{cor}[theorem]{Corollary}
\newtheorem{lem}[theorem]{Lemma}
\newtheorem{prop}[theorem]{Proposition}

\theoremstyle{definition}

\theoremstyle{remark}

\theoremstyle{remark}

\numberwithin{equation}{section}

\newcommand{\nc}{\newcommand}
\nc{\on}{\operatorname}
\nc{\ch}{\mbox{ch}}
\nc{\Z}{{\mathbb Z}}
\nc{\C}{{\mathbb C}}
\nc{\pone}{{\mathbb C}{\mathbb P}^1}
\nc{\pa}{\partial}
\nc{\F}{{\mathcal F}}
\nc{\arr}{\rightarrow}
\nc{\larr}{\longrightarrow}
\nc{\al}{\alpha}
\nc{\ri}{\rangle}
\nc{\lef}{\langle}
\nc{\W}{{\mathcal W}}
\nc{\la}{\lambda}
\nc{\ep}{\epsilon}

\nc{\su}{\widehat{{\mathfrak sl}}_2}
\nc{\sw}{{\mathfrak s}{\mathfrak l}}

\nc{\g}{{\mathfrak g}}
\nc{\h}{{\mathfrak h}}
\nc{\n}{{\mathfrak n}}
\nc{\N}{\widehat{\n}}
\nc{\G}{\widehat{\g}}
\nc{\De}{\Delta_+}
\nc{\gt}{\widetilde{\g}}
\nc{\Ga}{\Gamma}
\nc{\one}{{\mathbf 1}}
\nc{\z}{{\mathfrak Z}}
\nc{\zz}{{\mathcal Z}}
\nc{\Hh}{{\mathcal H}_\beta}
\nc{\qp}{q^{\frac{k}{2}}}
\nc{\qm}{q^{-\frac{k}{2}}}
\nc{\La}{\Lambda}
\nc{\wt}{\widetilde}
\nc{\qn}{\frac{[m]_q^2}{[2m]_q}}
\nc{\cri}{_{\on{cr}}}
\nc{\kk}{h^\vee}
\nc{\sun}{\widehat{\sw}_N}
\nc{\hh}{\widehat{\mathfrak h}}
\nc{\HH}{{\mathcal H}_{q,t}}
\nc{\ca}{\wt{{\mathcal A}}_{h,k}(\sw_2)}
\nc{\gl}{\widehat{{\mathfrak g}{\mathfrak l}}_2}
\nc{\el}{\ell}
\nc{\s}{{\mathbf s}}
\nc{\bi}{\bibitem}
\nc{\om}{\omega}
\nc{\WW}{\W_\beta}
\nc{\scr}{{\mathbf S}}
\nc{\ab}{{\mathbf a}}
\nc{\rr}{r}
\nc{\ol}{\overline}
\nc{\con}{qt^{-1} + q^{-1}t}
\nc{\den}{q^{\el-1} t^{-\el+1}+ q^{-\el+1} t^{\el-1}}
\nc{\ds}{\displaystyle}
\nc{\B}{B}
\nc{\A}{{\mathbb A}}
\nc{\GG}{{\mathcal G}}
\nc{\UU}{{\mathcal U}}
\nc{\MM}{{\mathcal M}}
\nc{\CC}{{\mathcal C}}
\nc{\GL}{{}^L G}
\nc{\dzz}{\frac{dz}{z}}
\nc{\Res}{\on{Res}}
\nc{\rep}{{\mathcal R}ep \;}
\nc{\uqg}{U_q \G}
\nc{\uqgg}{U_q \g}
\nc{\Fq}{{\mathbb F}_q}

\nc{\stimes}{\ltimes}
\nc{\K}{\hat{\mathcal K}}
\nc{\Ql}{\ol{\mathbb Q}_\ell}

\nc{\ga}{\gamma}
\nc{\PL}{{}^L P}
\nc{\E}{\mc E}
\nc{\mc}{\mathcal}
\nc{\mbf}{\mathbf}
\nc{\bb}{{\mathfrak b}}
\nc{\OO}{{\mc O}}
\nc{\Po}{{\mc P}}
\nc{\V}{{\mc V}}
\nc{\yy}{{\mc Y}}
\nc{\M}{\mathcal M}
\nc{\Coh}{{{\mathcal C}oh}}
\nc{\Cohn}{\Coh_n}
\nc{\f}{{\mathcal F}}
\nc{\si}{_E}
\nc{\Gaf}{{\mathbb G}_{a,\Fq}}
\nc{\KK}{{\mathfrak k}}

\nc{\PCr}{{ \bs P  (\C[x])^r   }}
\nc{\PCN}{{ \bs P  (\C[x])^N   }}
 
\nc{\sN}{sl_{2N+1}}

\nc{\Pzr}{{ \bs P(\C((x-z)))^r}}

\nc{\PzN}{{ \bs P(\C((x-z)))^N}}

\newcommand{\bean}{\begin{eqnarray}}
\newcommand{\eean}{\end{eqnarray}}
\newcommand{\be}{\begin{displaymath}}
\newcommand{\ee}{\end{displaymath}}
\newcommand{\bea}{\begin{eqnarray*}}   
\newcommand{\eea}{\end{eqnarray*}}
\newcommand{\bs}{\boldsymbol}
\newcommand{\Ref}[1]{{$($\ref{#1}$)$}}

\newcommand{\p}{\partial_x}

\newcommand{\Om}{\Omega}
\newcommand{\ox}{\otimes}

\def\Vb{L_\bullet}
\def\gl{\frak{gl}}
\def\glt{\gl_2}
\def\glN{\gl_N}
\def\glM{\gl_M}
\def\slN {\frak{sl}_N}
\def\V-{\>\hbox{$\=V\}$-}}

\newcommand{\NNN}{^{\langle N\rangle}}
\newcommand{\MMM}{^{\langle M\rangle}}
\newcommand{\Zp}{\Z_{\geq 0}}

\newcommand{\tl}{t^{\langle \bs n\rangle}}

\newcommand{\GR}{{\rm Gr}(X,N)}
\newcommand{\Wr}{{\rm Wr}}

\newcommand{\LMN}{\bs L_{\bs m}[\bs n]}
\newcommand{\LNM}{\bs L_{\bs n}[\bs m]}

\newcommand{\Sym}{{\rm Sym}}

\newcommand{\Hom}{{\rm Hom}}

\newcommand{\rank}{{\rm rank}}

\newcommand{\Gr}{{\rm Gr}}

\begin{document}

\title[Bispectral and $(\glN,\glM)$ Dualities]{Bispectral and $(\glN,\glM)$ Dualities}
\author[E. Mukhin, V. Tarasov, and A. Varchenko]
{E. Mukhin, V. Tarasov, and A. Varchenko}
\thanks{Research of E.M. was supported in part by NSF grant DMS-0140460.
Research of A.V. was supported in part by NSF grant DMS-0244579}

\address{E.M.: Department of Mathematical Sciences, Indiana University --
Purdue University Indianapolis, 402 North Blackford St, Indianapolis,
IN 46202-3216, USA,\newline mukhin@math.iupui.edu}
\address{V.T.: Department of Mathematical Sciences, Indiana University --
Purdue University Indianapolis, 402 North Blackford St, Indianapolis,
IN 46202-3216, USA,\newline vtarasov@math.iupui.edu, \ and\newline
\hglue\parindent St.~Petersburg Branch of Steklov Mathematical Institute,
Fontanka 27,\newline
St.~Petersburg, 191023, Russia, vt@pdmi.ras.ru}
\address{A.V.: Department of Mathematics, University of North Carolina
at Chapel Hill,\newline Chapel Hill, NC 27599-3250, USA, anv@email.unc.edu}


\maketitle

\begin{center} June, 2005
\end{center}

\begin{abstract}
Let $V = \langle\, p_{ij}(x)e^{\la_ix},\ i=1,\dots,n, \ j=1, \dots ,
N_i\, \rangle$ be a space of quasi-polynomials of dimension
$N=N_1+\dots+N_n$. Define the regularized fundamental operator of $V$
as the polynomial differential operator $D = \sum_{i=0}^N
A_{N-i}(x)\p^i$ annihilating $V$ and such that its leading coefficient $A_0$ is a
polynomial of the minimal possible degree. We construct a space of
quasi-polynomials $U = \langle \,q_{ab}(u)e^{z_au}\, \rangle$ whose
regularized fundamental operator is the differential operator
$\sum_{i=0}^N u^i A_{N-i}(\partial_u)$.  The space $U$ is constructed
from $V$ by a suitable integral transform. Our integral transform
corresponds to the bispectral involution on the space of rational
solutions (vanishing at infinity) to the KP hierarchy, see \cite{W}.

As a corollary of the properties of the integral transform
 we obtain a correspondence between
critical points of the two master functions associated with the $(\glN,\glM)$ dual
Gaudin models as well as between the corresponding Bethe vectors.

\end{abstract}

\section{Introduction}
Let $V = \langle\, p_{ij}(x)e^{\la_ix},\ i=1,\dots,n, \ j=1, \dots ,
N_i\, \rangle$ be a space of quasi-polynomials of dimension
$N=N_1+\dots+N_n$. Define the regularized fundamental operator of $V$
as the polynomial differential operator $D = \sum_{i=0}^N
A_{N-i}(x)\p^i$ annihilating $V$ and such that its leading coefficient $A_0$ is a
polynomial of the minimal possible degree. In this paper we construct a space of
quasi-polynomials $U = \langle \,q_{ab}(u)e^{z_au}\, \rangle$ whose
regularized fundamental operator is the differential operator
$\sum_{i=0}^N u^i A_{N-i}(\partial_u)$.  The space $U$ is constructed
from $V$ by a suitable integral transform, see Section \ref{sec int transform}.
Our integral transform
corresponds to the bispectral involution on the space of rational
solutions (vanishing to infinity) to the KP hierarchy, see \cite{W} and Section \ref{baker}.

As a corollary of the properties of the integral transform
we obtain a correspondence between
critical points of the two master functions associated with the $(\glN,\glM)$ dual
Gaudin models as well as between the corresponding Bethe vectors.

\medskip
\noindent
{\bf Example.}
Let $\bs n=(n_1,n_2)$ and $\bs m=(m_1,m_2)$ be two vectors of nonnegative integers such that
$n_1+n_2=m_1+m_2$. 
Consider two (master) functions,
\bea
&&
\Phi (t_1,\dots,t_{n_2}; \la_1,\la_2; z_1,z_2; \bs m) =
\exp\left(\la_1(m_1z_1 +m_2z_2) + (\la_2-\la_1)\sum_{j=1}^{n_2} t_j \right) \times
\\
&&
\phantom{aaaaaaaaaaa}
(z_1-z_2)^{m_1m_2}\ 
\prod_{j=1}^{n_{2}}
(t_j-z_1)^{-m_1}(t_j-z_2)^{-m_2}\!\!\!
\prod_{1\leq j<j'\leq n_2}\!\!\!
(t_j-t_{j'})^{2}\ ,
\\
&&
\Phi (s_1,\dots,s_{m_2}; z_1,z_2; \la_1,\la_2; \bs n) =
\exp\left(z_1(n_1\lambda_1 +n_2\lambda_2) + (z_2-z_1)\sum_{j=1}^{m_2} s_j \right) \times
\\
&&
\phantom{aaaaaaaaaaa}
(\la_1-\la_2)^{n_1n_2}\ 
\prod_{j=1}^{m_{2}}
(s_j-\la_1)^{-n_1}(s_j-\la_2)^{-n_2}\!\!\!
\prod_{1\leq j<j'\leq m_2}\!\!\!
(s_j-s_{j'})^{2}\ ,
\eea
depending on parameters $\la_1,\la_2,z_1,z_2$. Assume that the parameters are generic.
Both functions are symmetric with respect to permutations of their coordinates.

To the orbit of a 
critical point $(t_1,\dots, t_{n_2})$ of the first function assign the polynomial
$p_2(x)=\prod_{i=1}^{n_2}(x-t_i)$. Define the polynomial differential operator 
$D = A_0(x) \p^2 + A_1(x) \p + A_2(x)$ by the formula
\bea
D = (x-z_1)(x-z_2) 
(   \partial_x  - 
\ln'  ( \frac { e^{\la_1x}(x-z_1)^{m_1}(x-z_2)^{m_2} } { p_2(x) }  ) )
( \partial_x - \ln' (  p_{2}(x)e^{\la_2x} ) ) \ .
\eea
The differential equation $Df(x)=0$ has a solution $p_2(x)e^{\la_2x}$.
In Section \ref{sec m=n} we show that the differential equation
$ (u^2A_0(\partial_u) + u A_1(\partial_u)  + A_2(\partial_u))g(u)=0$ has a solution
of the form $e^{z_2u}\prod_{i=1}^{m_2}(u-s_i)$ and $(s_1,\dots,s_{m_2})$ is a 
critical point of the second function. This \linebreak
(bispectral) 
correspondence between the orbits of critical points of the two functions is reflexive.

To the orbit of a critical point $(t_1,\dots, t_{n_2})$ of the first function assign 
the  Bethe  vector
\bea
\omega (t_1, \dots , t_{n_2})\ =\!\!\! 
\sum_i
{C_i}\
 \frac{E_{21}^{n_2-i}v_{m_1}}{(n_2-i)!}\ \ox \frac{E_{21}^iv_{m_2}}{ i!}\ \ ,
\eea
where 
\bean\label{univ function}
C_i(t_1,\dots,t_{n_2})\ = \ \Sym\
 \prod_{j=1}^{n_2-i} 
\frac 1 {t_j - z_1}\
\prod_{j=1}^{i} 
\frac 1 {t_{n_2+j-i} - z_2}
\eean
and $\Sym \, f(t_1,\dots , t_{n})\ =\ \sum_{\sigma\in\Sigma_{n}} 
f(t_{\sigma(1)},\dots , t_{\sigma(n)})$, see \cite{FMTV}. 
This is a vector in the weight subspace of weight $[n_1,n_2]$ of
the tensor product $L_{m_1}\ox L_{m_2}$ of the 
irreducible $\glt$-modules with highest weights $(m_1,0)$ and $(m_2,0)$, respectively.
This vector is an eigenvector of the Gaudin Hamiltonians acting on the tensor product,
see Section \ref{KZ hamiltonians}.

Similarly, to the orbit of a critical point $(s_1,\dots, s_{m_2})$ of the second function 
assign the  Bethe  vector $\omega (s_1, \dots , s_{m_2})$ in 
the weight subspace of weight $[m_1,m_2]$ of the
tensor product $L_{n_1}\ox L_{n_2}$ of the 
irreducible $\glt$-modules with highest weights $(n_1,0)$ and $(n_2,0)$, respectively.

The Gaudin models on $(L_{m_1}\ox L_{m_2})[n_1,n_2]$ 
and $(L_{n_1}\ox L_{n_2})[m_1,m_2]$ are identified by the
$(\glt,\glt)$ duality of \cite{TV}.

In Section \ref{sec crit points} we show that the Bethe eigenvectors, corresponding to
the bispectral dual orbits of critical points, are identified by the $(\glt,\glt)$ duality
isomorphism.

\bigskip

The paper has the following structure. In Section \ref{spaces} we define the integral transform
of a space of quasi-polynomials and describe its properties. In Section \ref{special spaces}
we study special spaces of quasi-polynomials related to the $(\glN, \glM)$ duality.
In Section \ref{KZ hamiltonians} we describe the Gaudin model and sketch its Bethe
ansatz. In Section \ref{sec duality} we describe the $(\glN, \glM)$ duality and formulate
a conjecture about the bispectral correspondence of Bethe vectors under the
$(\glN, \glM)$ duality. In Section \ref{sec m=n} we prove the conjecture for $N=M=2$.
In Section \ref{baker} we indicate the relation of our integral transform to the 
bispectral involution on the space of rational solutions (vanishing at infinity)
to the KP hierarchy.
Section \ref{appendix} contains an appendix with an auxiliary result from 
representation theory.

\bigskip

The authors thank Yuri Berest for teaching them on the bispectral involution in the KP theory.

\section{ Spaces of Functions }\label{spaces}

\subsection{Fundamental operator}\label{subsection fund oper}
Let $\C$ be the complex line with coordinate $x$.
Let $V$ be a finite-dimensional complex vector space of meromorphic functions
on $\C$, dim $V = N$.

For $z\in \C$, define {\it the sequence of exponents of $V$ at $z$} as
the unique sequence of integers, $\bs e = \{e_1 < \dots < e_N\}$, with
the property: for $i=1,\dots,N$, there exists $f\in V$ such that $f$
has order $e_i$ at $z$.

We say that $z\in \C$ is a singular point of $V$ if the set of
exponents of $V$ at $z$ differs from the set $\{0,\dots, N-1\}$. We
assume that $V$ has at most finitely many singular points.

For meromorphic functions $f_1,\dots, f_i$ on $\C$, denote by
$\Wr(f_1,\dots,f_i)$ their Wronskian, the determinant of the $i\times
i$-matrix whose $j$-th row is $f_j, f_j^{(1)},...,f_j^{(i-1)}$.

Define the Wronskian of $V$, denoted by $\Wr_V$, as the Wronskian of a
basis of $V$.  The Wronskian of $V$ is determined up to multiplication
by a constant.

\bigskip

{\it The monic fundamental operator of $V$} is the unique monic linear
differential operator of order $N$ annihilating $V$. It is denoted by
$\bar D_V$.  We have \bean\label{fund oper} \bar D_V\ =\ \p^N + \bar
A_1\p^{N-1} + \dots + \bar A_N\ , \qquad \bar A_i\ =\ (-1)^i\, \frac
{\Wr_{V, i}}{\Wr_V}\ , \eean where $\p = d/dx$,\ {} $\Wr_V$ is the
Wronskian of a basis $f_1,\dots ,f_N \in V$, \ {} $\Wr_{V,i}$ is the
determinant of the $N\times N$-matrix whose $j$-th row is $f_j,
f_j^{(1)},...,f_j^{(N-i-1)},f_j^{(N-i+1)},\dots , f_j^{(N)}$.

Let $(z_1, \dots, z_m)$ be the subset of $\C$ of singular points of
$V$.  For $a=1,\dots , m$, let $M_a$ be the smallest natural number
such that all coefficients of the differential operator
$(x-z_a)^{M_a}\bar D_V$ have no poles at $x=z_a$. Then the
differential operator \bean\label{reg fund oper} D_V\ = \
\prod_{a=1}^m\ (x - z_a)^{M_a}\ \bar D_V \eean is called {\it the
regularized fundamental operator of $V$}. All coefficients of $D_V$
are regular in $\C$.

\begin{lem}\label{irreg of fund oper}
Let $z_a$ be a singular point of $V$ with exponents $\bs e$.  Let
$1\leq \bar M_a \leq N$ be the integer such that $e_i = i-1$ for
$i\leq N-\bar M_a$ and $e_{N-\bar M_a+1} \neq N-\bar M_a$. Then the
order of $\bar A_{\bar M_a}$ at $z_a$ is $-\bar M_a$ and the order of
$\bar A_i$ is not less than $-\bar M_a$ for $i > \bar M_a$.  \hfill
$\square$
\end{lem}
The proof follows from counting orders at $z_a$ of determinants $W_{V,i}$.

It is known that the order of any $A_i$ at $z_a$ is not less that $-i$.
In particular Lemma
\ref{irreg of fund oper}  implies that $\bar M_a =  M_a$.

\subsection{Example: Space of quasi-polynomials}\label{subs exponentials}
Let $N_1, \dots , N_n$ be natural numbers. Set $N = N_1 + \dots +
N_n$.  For $i = 1,\dots, n$,\ let \ $0 < n_{i1} < \dots < n_{iN_i}$\
be a sequence of positive integers.  For $i= 1,\dots , n,\ j = 1,
\dots , N_i$, \ let\ $p_{ij}\in \C[x]$\ be a polynomial of degree
$n_{ij}$.  Let $\lambda_1, \dots , \lambda_n$ be distinct complex
numbers.

Denote by $V$ the complex vector space spanned by functions \ $p_{ij}
e^{\lambda_ix}$,\ $i=1,\dots , n$,\linebreak $j=1,\dots , N_i$. The
dimension of $V$ is $N$.

Let $(z_1,\dots ,z_m)$ be the set of singular points of $V$. 
Let 
$$
\{\ 0 <  \dots < N - M_a -1 < N - M_a  +  m_{a 1} < \dots < N - M_a +  m_{a M_a}\ \}
$$ 
be the exponents of $V$ at $z_a$. Here $0<m_{a 1}<\dots < m_{a
M_a}$ and $M_a$ is an integer such that $1\leq M_a \leq N$.  Set $M =
M_1 + \dots + M_m$.

A space $V$ with such properties will be called {\it a space of
quasi-polynomials.}

\begin{lem}\label{Wr. identity}
We have
\bea
\sum_{a=1}^m\sum_{b=1}^{M_a}\ (\,m_{a b} +1-b\,) 
\ = \ \sum_{i=1}^n\sum_{j=1}^{N_i} \ (\,n_{i j}+1-j\,) \ .
\eea
\end{lem}
The lemma is proved by calculating the Wronskian of $V$.
\bigskip

Let $D_V = A_0 \p^N\!+A_1\p^{N-1}\!+\dots +A_N$
be the regularized fundamental operator of $V$. 

\begin{lem}\label{deg coeff}
${}$

\begin{enumerate}
\item[$\bullet$] We have
$
A_0 \ =\ \prod_{a=1}^m\ (x-z_a)^{M_a}.
$
\item[$\bullet$] All coefficients $A_i$ are polynomials in $x$ of
degree not greater than $M$.
\item[$\bullet$] Write $D_V = x^M B_0(\p) + x^{M-1} B_1(\p) + \dots + x B_{M-1}(\p)
+ B_M(\p)$ where $B_{i}(\p)$ is a polynomial in $\p$ with constant coefficients. Then
$
B_0\ =\ \prod_{i=1}^n \ (\p - \lambda_i)^{N_i}.
$
\item[$\bullet$] The polynomials $B_0(\p), \dots , B_M(\p)$ have no 
common factor of positive degree.
\end{enumerate}
\hfill
$\square$

\end{lem}

\subsection{Conjugate Space}\label{subsection conj space}
Let $V$ be a complex vector space of meromorphic functions
on $\C$ of dimension $N$ as in Section \ref{subsection fund oper}.

The complex vector space spanned by all functions of the form
$\Wr(f_1,\dots,f_{N-1})/\Wr_V$ with $f_i\in V$ has dimension $N$. It
is denoted by $V^\star$ and called {\it conjugate to} $V$.

The complex vector space spanned by all functions of the form
$f \prod_{a=1}^m(x-z_a)^{-M_a}$ with $f\in V^\star$
 has dimension $N$. It is  denoted by
$V^\dagger$ and  called {\it regularized conjugate to} $V$.

\begin{lem}\label{conj exponents} For $a=1, \dots , m$, if
 $\bs e$ are exponents  of $V$ at $z_a$, then
\bea
\bs e^\star = \{ - e_{N}-1+N <  - e_{N-1}-1+N < \dots <  - e_{1}-1+N\}
\eea
are exponents of $V^\star$ at $z_a$ and
\bea
\bs e^\dagger = \{ - e_{N}-1+N-M_a <  - e_{N-1}-1+N-M_a < \dots <  - e_{1}-1+N-M_a\}
\eea
are exponents of $V^\dagger$ at $z_a$.
\hfill
$\square$
\end{lem}

\bigskip

Let $D = \sum_i A_i \p^i$ be a differential operator with meromorphic coefficients.
The operator $D^* = \sum_i(-\p)^i A_i$  is called {\it
formal conjugate to} $D$.

\bigskip

\begin{lem}\label{conj mon oper}
Let $\bar D_V$ and $D_V$ be the monic and regularized fundamental
operators of $V$, respectively. Then $(\bar D_V)^*$ annihilates
$V^\star$ and $(D_V)^*$ annihilates $V^\dagger$.  \hfill $\square$
\end{lem}

\subsection{Integral transform}\label{sec int transform}
 Let $V$ be a space of quasi-polynomials as in Section \ref{subs
exponentials}.  For $a = 1, \dots , m$, let $\gamma_a$ be a small
circle around $z_a$ in $\C$ oriented counterclockwise. Denote by $U$
the complex vector space spanned by functions of the form
\bean\label{int transform} \hat f_a(u)\ =\ \int_{\gamma_a}\ e^{ux}\
f(x)\ dx , \eean where $a = 1, \dots , m$, \ $f \in V^\dagger$. The
vector space $U$ is called {\it bispectral dual to $V$}.

\begin{theorem}\label{thm int transform}
${}$

\begin{enumerate} 
\item[(i)] The space $U$ has dimension $M$.
\item[(ii)] For $a=1,\dots, m$,\ $b=1, \dots , M_a$, the space $U$ contains a function
of the form $q_{a b}e^{z_a u}$ where $q_{a b}$ is a polynomial in $u$ of degree $m_{a b}$.
\item[(iii)] Let $D_V = \sum_{i=1}^M\sum_{j=1}^N A_{i j}x^i\partial_x^j$ 
be the regularized fundamental operator of $V$ where $A_{i j}$ are suitable complex numbers.
Then
\bean\label{conj reg fund oper}
\sum_{i=1}^M\sum_{j=1}^N\ A_{i j}\ u^j \partial_u^i
\eean
 is the regularized fundamental operator of $U$.
\item[(iv)]
The set $( \lambda_1, \dots , \lambda_n )$ 
is the subset in  $\C$  of singular points of $U$.
\item[(v)]
For $i = 1,\dots , n$,\ the 
 set 
\bea
\{0 < \dots < M-N_i-1 < M-N_i + n_{i 1}< M-N_i+n_{i 2}< \dots < M-N_i+n_{i N_i}  \}
\eea
is the set of exponents of
$U$ at $\lambda_i$.
\end{enumerate}
\end{theorem}

\begin{cor}\label{cor involutivity} 
The space $V$ is bispectral dual to $U$.

\end{cor}

{\it Proof of the theorem.}
The exponents of $V^\dagger$ at $z_a$ are 
$$
\{-m_{a,M_a}-1 < \dots
< -m_{a1}-1< 0< 1<\dots < N - M_a -1 \}\ .
$$ 
Integral \Ref{int transform} is nonzero only if $f$ has a pole
at $x=z_a$. If $f$ has a pole at $x=z_a$ of order $-m_{ab}-1$, then
integral \Ref{int transform} has the form $q_{ab} e^{z_a u}$ where 
$q_{ab}$ is a polynomial in $u$ of degree $m_{ab}$. This proves (i) and (ii).

It is clear that the operator
$D^\dagger = \sum_{i=1}^M\sum_{j=1}^N\ A_{i j}\ u^j \partial_u^i$\
annihilates $U$. Write 
\bea
D^\dagger\ =\ B_0(u)\partial_u^M + \dots + B_{M-1}(u) \partial_u + B_M(u)
\eea
where $B_a(u)$ are polynomials in $u$ with constant coefficients. 
Lemma \ref{deg coeff} implies that 
$B_0(u)= \prod_{j=1}^{n}(u-\lambda_j)^{N_j}$ and
the polynomials $B_0, \dots , B_M$ have no common factor of positive degree. 
Therefore, $D^\dagger$ is the
 regularized fundamental operator of $U$. Part (iii) is proved.


Using Lemma \ref{irreg of fund oper}  we get $N_i \leq M$ for all $i$.

The Wronskian of $V$ has the form  $\Wr_V = p\,e^{x\sum_{i=1}^n N_i \lambda_i}$ 
where $p$ is a polynomial in $x$ of degree 
$ \sum_{i=1}^n\sum_{j=1}^{N_i} \ (\,n_{i j}+1-j\,)$.

The functions $p_{ij}e^{\lambda_i x}$, $i=1,\dots,n,\
j=1, \dots , N_i$,\ form a basis of $V$. For given $i, j$, order
all the basis functions of $V$ except the function
$p_{ij}e^{\lambda_i x}$. Denote by $\Wr_{ij}$ the Wronskian 
of this ordered set of functions. The functions
$\Wr_{ij}/\Wr_V\prod_{a=1}^m(x-z_a)^{M_a}$ form a basis in $V^\dagger$.
The function $\Wr_{ij}/\Wr_V\prod_{a=1}^m(x-z_a)^{M_a}$ has the form
$r_{ij}e^{-\lambda_i x}$ where $r_{ij}$ is a rational function in $x$.
The function $r_{ij}$ has zero at $x=\infty$ of  order $M + n_{ij}- N_i + 1$.
Notice that this integer is positive.

Let $\gamma_\infty$ be a big circle around infinity oriented clockwise.
The value of the function
\bea
\sum_{a=1}^m \int_{\gamma_a}\ e^{ux}\ e^{-\lambda_i x} \ r_{ij}(x)  \ dx \
=\
\int_{\gamma_\infty}\ e^{ux}\ e^{-\lambda_i x}\ r_{ij}(x)  \ dx 
\eea
and its first $M + n_{ij} - N_i - 1$ derivatives at $u=\lambda_i$ is zero,
but the value at $u=\lambda_i$ of its
$(M-N_i + n_{ij})$-th derivative is not zero.  This remark together with Lemma
\ref{Wr. identity} proves (iv) and (v). 
\hfill
$\square$

\section{Special spaces}\label{special spaces}

\subsection{Spaces of a $(\bs \la, \bs z, \bs n, \bs m)$-type}
Let $N > 1$ be a natural number.
Let $\bs n = \linebreak
(n_1, \dots , n_N)$ be a vector of nonnegative integers. Let
$\bs \lambda = (\lambda_1, \dots , \lambda_N)$ be a vector with distinct complex coordinates.
For $i=1,\dots , N$, let $p_i\in \C[x]$ be a polynomial of degree $n_i$.
Denote by $V$ the complex vector space spanned by functions \ $p_{i} e^{\lambda_ix}$,\
$i=1,\dots , N$. The dimension of $V$ is $N$.

Let $\bs z = (z_1,\dots ,z_M)$, $M>1$,\ be a  subset in $\C$ containing all singular points of $V$. 
Assume that for $a = 1,\dots , M$, the set of exponents of $V$ at
$z_a$ has the form 
$$
\{ 0 < 1 < \dots < N-2< N-1+m_{a}\}\ ,
\qquad
m_a \geq 0\ .
$$
We have 
$$
\sum_{i=1}^N \ n_i\ =\ \sum_{a=1}^M\ m_a\ .
$$
We call the pair $(V, \bs z)$ {\it a space of the $(\bs \la, \bs z, \bs n, \bs m)$-type} or
{\it a special space}.

Let  $\bar D_V$ be the monic fundamental operator of $V$. The operator
\bea
\tilde D_V\ = \ \prod_{a=1}^M\ (x - z_a)\ \bar D_V 
\eea
is called {\it the special fundamental operator } of the special space $(V,\bs z)$.
By Lemma \ref{deg coeff}, all coefficients of the differential operator
$\tilde D_V$ are polynomials in $x$ of
degree not greater than $M$.  If we write 
$$
\tilde D_V = x^M B_0(\p) + x^{M-1}
B_1(\p) + \dots + x B_{M-1}(\p) + B_M(\p),
$$
where $B_i$ is a polynomial in $\p$ with constant coefficients, 
 then $B_0 = \prod_{i=1}^N (\p - \lambda_i)$.

\subsection{Special integral transform}
Let $(V,\bs z)$ be a space of a $(\bs \la, \bs z, \bs n, \bs m)$-type.
For $a = 1, \dots , M$, let $\gamma_a$ be a small circle around $z_a$ in $\C$
 oriented counterclockwise. Denote by $U$ the complex vector space
spanned by functions of the form
\bean\label{int transform}
\hat f_a(u)\ =\ \int_{\gamma_a}\ e^{ux}\ f(x)\ \prod_{a=1}^M(x-z_a)^{-1}\ dx\ ,
\eean
where $a = 1, \dots , m$, \ $f \in V^\star$. 

\begin{theorem}\label{thm special int transform}
${}$

\begin{enumerate} 
\item[(i)] The space $U$ has dimension $M$.
\item[(ii)] For $a=1,\dots, M$, the space $U$ contains a function
of the form $q_{a}e^{z_a u}$ where $q_{a}$ is a polynomial in $u$ of degree $m_{a}$.
\item[(iii)] The subset $\bs\la = (\la_1, \dots , \la_N)$ of $\C$ contains all singular
points of $U$.
\item[(iv)]
For $i = 1,\dots , N$,\ the 
 set 
\bea\{0 < \dots < M-2 < M - 1 + n_{i}  \}
\eea
is the set of exponents of
$U$ at $\lambda_i$. Thus $(U, \bs \la)$ is a space of
the  $(\bs z, \bs \la, \bs m, \bs n)$-type.
\item[(v)] Let $\tilde D_V = \sum_{i=1}^M\sum_{j=1}^N A_{i j}x^i\partial_x^j$ 
be the special fundamental operator of $(V, \bs z)$ where $A_{i j}$ are suitable complex numbers.
Then
\bean\label{conj reg fund oper}
\sum_{i=1}^M\sum_{j=1}^N\ A_{i j}\ u^j \partial_u^i
\eean
 is the special fundamental operator of $(U, \bs \la)$.
\end{enumerate}
\end{theorem}

The proof is similar to the proof of Theorem \ref{thm int transform}.

The special space $(U,\bs\la)$ is called {\it  special bispectral dual} to $(V,\bs z)$.

\begin{cor}\label{cor involutivity} 
The space $(V, \bs z)$ is special bispectral dual to $(U, \bs \la)$.

\end{cor}

\subsection{Special spaces and critical points}\label{sec master function}
Let $(V, \bs z)$ be a space of a $(\bs \la, \bs z, \bs n, \bs m)$-type. 
We construct the associated master function
as follows.
Set
$$
\bar n_i\ =\ n_{i+1} +\dots + n_N\ ,
\qquad
i=1,\dots, N-1\ .
$$
Consider new
$\bar n_1+\dots +\bar n_{N-1}$  auxiliary variables
$$ 
t^{\langle \bs n\rangle} \ = \ (t^{(1)}_1, \dots , t^{(1)}_{\bar n_1},\ 
t^{(2)}_1, \dots , t^{(2)}_{\bar n_2},\ 
\dots ,\ t^{(N-1)}_1,
\dots , t^{(N-1)}_{\bar n_{N-1}})\ .
$$
Define {\it the master function}
\bean
&&
\Phi (t^{\langle \bs n\rangle};\bs\la; \bs z; \bs m) =
\exp\left(\la_1\sum_{a=1}^Mm_az_a +
\sum_{i=1}^{N-1}(\la_{i+1}-\la_i)\sum_{j=1}^{\bar n_i} t^{(i)}_j\right)\!
\prod_{1\leq a<b\leq M}(z_a-z_b)^{m_am_b}
\notag
\\
\label{master function}
&&
\phantom{aa}
\times 
\prod_{a=1}^M
\prod_{j=1}^{\bar n_{1}}
(t^{(1)}_j-z_a)^{-m_a}\
\prod_{i=1}^{N-1}
\prod_{j<j'}
(t^{(i)}_j-t^{(i)}_{j'})^{2}
\prod_{i=1}^{N-2}
\prod_{j=1}^{\bar n_{i}}
\prod_{j'=1}^{\bar n_{i+1}}
(t^{(i)}_j-t^{(i+1)}_{j'})^{-1}\ .
\eean
The master function is symmetric with respect to the group
$\Sigma_{\bs{ \bar n}}=\Sigma_{\bar  n_1}\times \dots\times \Sigma_{\bar  n_{N-1}}$
of permutations of variables $t^{(i)}_j$ preserving the  upper index.

The $\Sigma_{\bs{ \bar n}}$-orbit of a point $t^{\langle \bs n\rangle}
 \in \C^{\bar n_1+\dots+\bar n_{N-1}}$ is uniquely determined by the
 $N-1$-tuple $\bs y = (y_1, \dots , y_{N-1})$ of polynomials in $x$,\ 
 where 
\bea 
y_i\ =\ \prod_{j=1}^{\bar n_{i}}\ (x - t^{(i)}_j)\ ,
 \qquad i=1,\dots, N-1\ .  
\eea 
We say that $\bs y$ {\it represents}
the orbit. Each polynomial of the tuple is
 considered up to multiplication by a number since we are interested in the
 roots of the polynomial only.

\medskip

We say that  $t^{\langle \bs n\rangle} \in \C^{\bar n_1+\dots+\bar n_{N-1}}$
is  {\it admissible} if the value  $\Phi (t^{\langle \bs n\rangle};\bs \la;\bs z; \bs m)$
is well defined and is not zero.

A point $t^{\langle \bs n\rangle}$ is admissible if and only if the associated tuple has
the following properties. The polynomial $y_1$ has no roots in $( z_1,\dots , z_M )$
and for all $i$,\
the polynomial $y_i$ has no multiple roots and no common roots with $y_{i-1}$ or $y_{i+1}$.
Such tuples are called admissible.

\medskip

The space $V = \langle p_1 e^{\lambda_1x}, \dots , p_N e^{\lambda_Nx} \rangle$
determines the $N-1$-tuple
\ $\bs y^{V} = (y^V_1, \dots , y^V_{N-1})$ of polynomials in $x$,\ {} where
\bea
y^V_i\ =\ e^{-(\lambda_{i+1}+\dots +\lambda_N)x}\ \Wr (p_{i+1}e^{\lambda_{i+1}x}, 
\dots , p_Ne^{\lambda_Nx})\ ,
\qquad
i=1,\dots, N-1\ .
\eea
We call the special space $(V, \bs z)$ {\it admissible} if the tuple $\bs y^V$ is admissible.

\begin{theorem}\cite{MV1, MV3}\label{space - crit pt}
${}$

\begin{enumerate}
\item[(i)]
Assume that the special space $(V, \bs z)$ is admissible. 
Then the tuple $\bs y^V$ represents the orbit of
a critical point of the master function.
\item[(ii)]
Assume that $t^{\langle \bs n \rangle}$ is admissible and $t^{\langle
\bs n\rangle}$ is a critical point of the master function. Let
$\bs y=(y_1,\dots,y_{N-1})$ be the tuple representing the orbit of $\tl$. 
Then the differential operator
\bea
\bar D & = & 
(   \partial_x  - 
\ln'  ( \frac { e^{\la_1x}\prod_{a=1}^M (x-z_a)^{m_a} } { y_1 }  ) )
( \partial_x - \ln' ( \frac {y_1e^{\la_2x} }{y_2}) )
( \partial_x - \ln' ( \frac {y_2e^{\la_3x} }{y_3}) )
\\
& \dots &
 ( \partial_x  -  \ln'  ( \frac { y_{N-2}e^{\la_{N-1}x} }{ y_{N-1} }  ) ) \
( \partial_x - \ln' (  y_{N-1}e^{\la_Nx} ) ) 
\eea
is the monic fundamental differential operator of a special space $(V, \bs z)$ of
the $(\bs \la, \bs z, \bs n, \bs m)$-type.
\item[(iii)] The correspondence  between 
admissible spaces of the 
$(\bs \la, \bs z, \bs n, \bs m)$-type and orbits 
of admissible critical points of the master function  described in parts 
{\rm (i), (ii)} is reflexive.
\end{enumerate}
\end{theorem}

This theorem establishes a one-to-one correspondence between admissible
spaces of the $(\bs \la, \bs z, \bs n, \bs m)$-type and orbits of admissible critical points
of the master function.

\subsection{Finiteness of admissible critical points}

\begin{lem}\label{finiteness}
For generic $\bs \la$ the master function has only finitely many admissible critical points.
\end{lem}
\begin{proof}
The admissible critical point equations are 
\bea
\sum_{j'=1, \ {} j' \ne j}^{\bar {n}_1}
\frac {2}{t^{(1)}_j-t^{(1)}_{j'}}
-
\sum_{j'=1}^{\bar n_{2}}
\frac 1{t^{(1)}_j-t^{(2)}_{j'}}
-
\sum_{a=1}^M \frac{m_a}{t^{(1)}_j-z_a} = \la_1-\la_2\ ,
&&
j=1, \dots , \bar n_1\ ,
\\
\sum_{j'=1, \ {} j' \ne j}^{\bar {n}_i}
\frac {2}{t^{(i)}_j-t^{(i)}_{j'}}
-
\sum_{j'=1}^{\bar n_{i-1}}
\frac 1{t^{(i)}_j-t^{(i-1)}_{j'}}
-
\sum_{j'=1}^{\bar n_{i+1}}
\frac 1{t^{(i)}_j-t^{(i+1)}_{j'}}  = \la_i-\la_{i+1}\ ,
&&
i=1,\dots, N-1,\ \ 
\\
&&
j=1, \dots , \bar n_i\ .
\eea
Assume that the admissible critical 
set contains an algebraic curve. Then   some of coordinates
$t^{(i)}_j$ tend to infinity along the curve. Add all equations corresponding to
those coordinates and take the limit. Then the left hand side of the resulting equation 
is zero while the right hand side is  $\sum_i c_i(\lambda_i-\la_{i+1})$
where $c_i$ is the number of coordinates $t^{(i)}_j$ which tend to infinity along the curve.
If the vector
$\bs \la$ is such that the numbers $\sum_i c_i(\lambda_i-\la_{i+1})$ are all nonzero,
then all admissible critical points are isolated.
\end{proof}

\subsection{The number of orbits of admissible critical points}
\label{number}
Set $\bar n = n_1+\dots + n_N$. Consider the complex vector space $X$ spanned by
functions $x^je^{\la_ix}$, \ $i=1, \dots , N$, $j = 0,\dots ,n_i$. The space $X$ is of dimension
$\bar n + N$. 

\begin{lem} The Wronskian  of $X$ is $e^{(n_1\lambda_1 + \dots + n_N\la_N)x}$.
\hfill
$\square$
\end{lem}

For $z\in\C$, introduce a complete flag $\bs F(z)$ in $X$,
\bea
\bs F(z)\ =\ \{ 0=F_0(z) \subset F_1(z) \subset
\dots \subset F_{\bar n+N}(z)=X \}\ ,
\eea
 $F_k(z)$ consists of all $f\in X$ which have zero at $z$ of order not less that
$\bar n + N -k$. By the lemma, $F_i(z)$ has dimension $k$.

Define two complete flags of $X$ {\it at infinity}. 

Say that $x^je^{\la_ix}\ <_1\
x^{j'}e^{\la_{i'}x}$ if $i<i'$ or $i=i'$ and $j<j'$. Set
\bea
\bs F(\infty_1)\ =\ \{ 0=F_0(\infty_1) \subset F_1(\infty_1) \subset
\dots  \subset F_{\bar n+N}(\infty_1)=X \}\ ,
\eea
where $F_k(\infty_1)$ is spanned by $k$ smallest elements with respect to $<_1$.

Say that $x^je^{\la_ix}\ <_2\
x^{j'}e^{\la_{i'}x}$ if $i>i'$ or $i=i'$ and $j<j'$. Set
\bea
\bs F(\infty_2)\ =\ \{ 0=F_0(\infty_2) \subset F_1(\infty_2) \subset
\dots  \subset F_{\bar n+N}(\infty_2)=X \}\ ,
\eea
where $F_k(\infty_2)$ is spanned by $k$ smallest elements with respect to $<_2$.

\bigskip

Denote by $\GR$  the Grassmannian manifold of $N$-dimensional vector subspaces of $X$.
Let $\bs F$ be a complete flag of $X$,
$$
\bs F\   = \ \{0=F_0
\ \subset F_1\ \subset \ \dots\ \subset \ F_{\bar n +N}\ =\ X \}\ .
$$ 
A {\it ramification sequence} is a sequence 
$(c_1,\dots,c_N)\in \Z^N$ such that $\bar n \geq c_1 \geq \dots \geq c_{N} \geq
0.$ For a ramification sequence $\bs c =(c_1,\dots,c_{N})$
 define the Schubert cell
\begin{align}
&
\Omega^o_{\bs c}(\bs F) \ =\ \{V\in \GR \mid \dim(V\cap F_u) = \ell ,\
\notag
\\
&
\phantom{aaaaaaaaaa}
\bar n +\ell-{c}_{\ell} \leq u < \bar n + \ell + 1 - {c}_{\ell+1},
\ \text{}\ell = 0, \dots,N \}\ ,
\notag
\end{align}
where ${c}_0 = \bar n$,\ {}\ ${c}_{N+1} = 0$. The cell
$\Omega^o_{\bs c}(\bs F)$ is a smooth connected variety. 
The closure of $\Omega^o_{\bs c}(\bs F)$ is denoted by 
$\Omega_{\bs c}(\bs F)$. The codimension 
of $\Omega^o_{\bs c}(\bs F)$ is
$$
|\bs c|\ =\ {c}_1\ +\ 
{c}_2\ +\ \dots\ +\ {c}_{N}\ . 
$$ 
Every 
$N$-dimensional vector subspace of $X$  belongs to a unique Schubert cell
$\Omega^o_{\bs c}(\bs F)$. 

\bigskip

For $a=1,\dots , M$, define the ramification sequence 
$\bs c(a)= (m_a,0,\dots,0)$. 
Define the ramification sequences
\bea
\bs c(\infty_1)& =& (n_2+\dots+n_N, n_3+\dots+ n_N,\dots,n_{N},0)\ ,
\\
\bs c(\infty_2) & = & (n_1+\dots+n_{N-1}, n_1+\dots+ n_{N-2},\dots,n_{1},0)\ .
\eea

\begin{lem} ${}$

\begin{enumerate}
\item[$\bullet$] We have
\bea
&&
\sum_{a=1}^M \ {\rm codim}\ \Omega^o_{\bs c(a)}(\bs F(z_a))\ +\
{\rm codim}\ \Omega^o_{\bs c(\infty_1)}(\bs F(\infty_1)) 
 +
\\ 
&&
\phantom{aaaaaaaaaaaaaa}
{\rm codim}\ \Omega^o_{\bs c(\infty_2)}(\bs F(\infty_2)) 
= 
\dim\ \GR\ =\ N \bar n\ .
\eea
\item[$\bullet$]
Let $V\in\GR$. The pair $(V, \bs z)$ is a space of the $(\bs \la, \bs z, \bs n, \bs m)$-type, if and only if
$V$ belongs to the intersection of $M+2$ Schubert cells
\bea
&&
\Omega^o_{\bs c(1)}(\bs F(z_1)) \cap \Omega^o_{\bs c(2)}(\bs F(z_2)) \cap
\dots \phantom{aaaaaaaaaaaaaa}
\\
&&
\phantom{aaaaaaaa}
\cap
\Omega^o_{\bs c(M)}(\bs F(z_M)) \cap
\Omega^o_{\bs c(\infty_1)}(\bs F(\infty_1)) \cap
\Omega^o_{\bs c(\infty_2)}(\bs F(\infty_2)) \ .
\eea
\end{enumerate}
\hfill
$\square$
\end{lem}

According to Schubert calculus,
the multiplicity of the intersection of Schubert cycles
\bean\label{intersection}
&&
\Omega_{\bs c(1)}(\bs F(z_1)) \cap \Omega_{\bs c(2)}(\bs F(z_2)) \cap
\dots 
\notag
\\
&&
\phantom{aaaaaa}
\cap
\Omega_{\bs c(M)}(\bs F(z_M)) \cap
\Omega_{\bs c(\infty_1)}(\bs F(\infty_1)) \cap
\Omega_{\bs c(\infty_2)}(\bs F(\infty_2)) \ 
\eean
can be expressed in representation-theoretic terms as follows.

For a ramification sequence $\bs c$ denote by $L_{\bs c}$ the finite
dimensional irreducible \linebreak
$\glN$-module with highest weight $\bs c$.
Any $\glN$-module $L$ has a natural structure of an $\slN$-module
denoted by $\widetilde L$.  By \cite{Fu},
the multiplicity of the intersection in
\Ref{intersection} is equal to the multiplicity of the trivial
$\slN$-module 
in the tensor product of $\slN$-modules
\bean\label{slN intersection}
\widetilde L_{\bs c(1)} \ox\dots \ox\widetilde L_{\bs c(M)} \ox
\widetilde L_{\bs c(\infty_1)} \ox \widetilde L_{\bs c(\infty_2)} \ .
\eean

\begin{prop}\label{lem on multiplicity}
The multiplicity of the trivial
$\slN$-module in the tensor product \Ref{slN intersection} is equal to the dimension 
of the weight subspace of weight $[n_1, \dots , n_N]$ in the tensor product of $\glN$-modules
\bean\label{glN intersection}
 L_{\bs c(1)} \ox\dots  \ox L_{\bs c(M)}\ .
\eean
\end{prop}

The proposition is proved in the appendix.

\begin{cor}\label{cor bound}
For generic $\bs \la$, the number of orbits of
admissible critical points of
the master function  $\Phi (t^{\langle \bs n\rangle};\bs\la; \bs z; \bs m)$
is not greater than
the dimension of the weight space
\linebreak
$ (L_{\bs c(1)} \ox\dots \ox L_{\bs c(M)})[n_1,\dots,n_N]$.
\end{cor}

\noindent
{\bf Remark.} Similarly to \cite{BMV}, one can show that the multiplicity
of an admissible critical point $t^{\langle \bs n \rangle}$ of the master function
is equal to the multiplicity of the corresponding point in the intersection
\Ref{intersection}. Thus the number of orbits of admissible critical points 
counted with multiplicities is not greater than the dimension of 
$ (L_{\bs c(1)} \ox\dots \ox L_{\bs c(M)})[n_1,\dots,n_N]$.

\section{KZ and dynamical Hamiltonians and critical points}\label{KZ hamiltonians}

\subsection{KZ and dynamical Hamiltonians}
Let $E_{ij}$, $i,j=1,\dots , N$, be the standard generators of the 
complex Lie algebra $\frak {gl}_N$.

We have the root decomposition \,$\glN=\n^+\oplus\h\oplus\n^-$ where
$$
\n^+\>=\,\oplus _{i<j}\ \C\cdot E_{i j}\,,\qquad
\h\,=\,\oplus_{i=1}^N \ \C \cdot E_{i i}\,,\qquad
\n^-\>=\,\oplus _{i>j}\ \C \cdot E_{i j}\,.
$$
The Casimir element is the element
$\Omega\,=\, \sum_{i,j = 1}^N E_{i j}\ox E_{j i}\ \in\glN^{\ox 2}$.

Let $Y = Y_1\ox\dots\ox Y_M$ be the tensor product of
$\glN$-modules.

The {\it KZ Hamiltonians}
$H_a(\bs\la,\bs z)$, $a=1,\dots, M$, acting on $Y$-valued
functions of \ $\la=(\la_1,\dots , \la_N), \ \bs z=(z_1,\dots, z_M) $ \  are defined
by the formula
$$
H_a(\bs\la, \bs z)\,=\,
\sum_{b=1,\ b \ne a}^M
\frac{\Om^{(a b)}}{z_a - z_b} \ +\ \sum_{i=1}^N\,\la_i\, E_{i i}^{(a)}\ .
$$
Here $\Om^{(a b)} : Y \to Y$ acts as $\Om$ on $Y_a\ox Y_b$,
and as the identity on other tensor factors. Similarly, 
$E_{i i}^{(a)}$ acts as $E_{i i}$ on $Y_a$ and as the identity on other
factors.

The {\it dynamical Hamiltonians} $G_i(\bs\la, \bs z)$,
$i=1,\dots,N$, acting on $Y$-valued functions of $\bs \la, \bs z$  are defined by the formula
\bea
G_i(\bs \la,\bs z)\,  =\
\sum_{j=1,\ j \ne i}^N
\frac {E_{i j} E_{j i}- E_{i i}}{\la_i -\la_j}\ 
+ \ \sum_{a=1}^M\ z_a\, E_{i i}^{(a)}\ .
\eea
The KZ and dynamical Hamiltonians commute  \cite{FMTV},
$$
[H_a(\bs\la,\bs z) , H_b(\bs \la, \bs z)]\ =\ 0\ ,
\qquad [H_a(\bs\la,\bs z)\>,G_i(\bs\la,\bs z)]\ =\ 0\ ,
\qquad
[G_i(\bs \la,\bs z) ,G_j(\bs\la,\bs z)]\ =\ 0\ ,
$$
for $a, b = 1,\dots, M$, and $i, j = 1,\dots,N$.

{\it The Gaudin diagonalization problem} is to diagonalize simultaneously
the KZ Hamiltonians $H_a$, $a=1,\dots, M$,\
and the dynamical Hamiltonians  $ G_i$, 
$i=1,\dots, N$, for given $\bs\la, \bs z$. 
The Hamiltonians preserve the weight decomposition of $Y$ 
and the diagonalization problem can be considered on a given weight subspace of $Y$.

\subsection{Diagonalization and critical points}
For a nonnegative integer $m$, denote by $L_m$ the irreducible  
$\glN$-module with highest weight $(m,0,\dots,0)$. 

Let $\bs n=(n_1,\dots,n_N)$ and $\bs m=(m_1,\dots,m_M)$ be vectors of nonnegative integers
with $\sum_{i=1}^N n_i = \sum_{a=1}^M m_a$.

Consider the tensor product $L_{m_1}\ox\dots\ox L_{m_M}$ and its weight subspace
\bea
\LMN\ =\ (L_{m_1}\ox\dots\ox L_{m_M})[n_1,\dots,n_N]\ .
\eea
Let $\bs \la \in \C^N$, $\bs z\in\C^M$. Assume that each of $\bs\la$ and 
$ \bs z$ has distinct coordinates.
Consider the Hamiltonians $H_a(\bs \la,\bs z)$, 
$G_i(\bs \la,\bs z)$ 
acting on $\LMN$.
{\it The Bethe ansatz method} is a method
to construct common eigenvectors of the Hamiltonians.

As in Section \ref{sec master function} consider the space
$\C^{\bar n_1+\dots +\bar n_{N-1}}$ with coordinates
$t^{\langle \bs n\rangle}$. Let $\Phi(\,\cdot\,;\bs\la;\bs z;\bs m)$ be the master function on
$\C^{\bar n_1+\dots +\bar n_{N-1}}$ defined in \Ref{master function}.

In Section 3 of \cite{FMTV}, a certain  $\LMN$-valued 
rational function 
$\omega : \C^{\bar n_1+\dots +\bar n_{N-1}}\to \LMN$\
is constructed. It is called {\it the universal rational function}.
Formulas for that function see also in \cite{RSV}.
For $N=2$ formulas for the universal rational 
function see in Section \ref{proof of conj}.

The universal rational function is well defined for admissible
$t^{\langle \bs n\rangle}$ \cite{RSV}. The universal 
rational function is symmetric with respect to
the  $\Sigma_{\bs{ \bar n}}$-action. 

\begin{theorem} \cite{RV, FMTV}\label{thm RV} If
$t^{\langle \bs n\rangle}$ is an admissible nondegenerate
critical point of the master function, then 
$\omega(t^{\langle \bs n\rangle})\in \LMN$ is an eigenvector of the Hamiltonians
with eigenvalues given by the derivatives of the logarithm of the master 
function with respect to the corresponding parameters
$z_a$ or $\la_i$,
\bean\label{eigenvalue}
H_a(\bs \la,\bs z)\ \omega(t^{\langle \bs n\rangle})
& = &
\frac{\partial }{\partial z_a} {\rm log} \Phi (t^{\langle \bs n\rangle};\bs \la; \bs z;\bs m)\
\omega(t^{\langle \bs n\rangle})\ ,
\qquad
a=1,\dots , M\ ,
\notag
\\
G_i(\bs \la,\bs z) \ \omega(t^{\langle \bs n\rangle})
& = &
\frac{\partial }{\partial \la_i} {\rm log} \Phi
(t^{\langle \bs n\rangle};\bs \la; \bs z;\bs m)\
 \omega(t^{\langle \bs n\rangle})\ ,
\qquad
i=1,\dots , N\ .
\eean

\end{theorem}

See comments on this theorem in \cite{V}.

The eigenvector $\omega(t^{\langle \bs n\rangle})$ is called {\it a Bethe eigenvector}.

By Corollary \ref{cor bound}, the number of Bethe eigenvectors is not greater than the dimension
of the space $\LMN$. The Bethe ansatz conjecture says that all eigenvectors of the Hamiltonians
are the Bethe eigenvectors.
Counterexamples to the Bethe ansatz conjecture see in \cite{MV2}.

\bigskip

According to previous discussions, eigenvectors of Hamiltonians on $\LMN$
are related to critical points of the master function $
\Phi( t^{\langle \bs n\rangle}; \bs \la; \bs z;\bs m)$. The critical
points of $ \Phi(t^{\langle \bs n\rangle};\bs \la; \bs z;\bs m)$ are
related to spaces of the $(\bs \la, \bs z, \bs n, \bs m)$-type. The
spaces of the \linebreak
$(\bs \la, \bs z, \bs n, \bs m)$-type are special bispectral
dual to spaces of the $(\bs z, \bs \la, \bs m, \bs n)$-type. The
spaces of the $(\bs z, \bs \la, \bs m, \bs n)$-type are related to
critical points of the master function $ \Phi(t^{\langle \bs
m\rangle}; \bs z; \bs \la; \bs n)$ which in its turn are related to
eigenvectors of the Hamiltonians acting on $\LNM$. As as a result of
this chain of relations, the eigenvectors of the Hamiltonians on
$\LMN$ and $\LNM$ must be related.

Indeed this relation is given by the $(\glN\>,\glM)$ duality.

\section{The $(\glN\>,\glM)$ duality for KZ and dynamical Hamiltonians}
\label{sec duality}

\subsection{The $(\glN\>,\glM)$ duality}\label{sec duality}
The Lie algebra $\glN$ acts on  $\C[x_1,\dots, x_N]$ by differential operators
$\displaystyle E_{ij}\mapsto x_i\frac \partial{\partial x_j}$. Denote this 
$\glN$-module by $\Vb$.
 Then
$$
\Vb\,=\,\bigoplus_{m=0}^\infty\,L_{m}\,,
$$
the submodule $L_{m}$ being spanned by homogeneous polynomials of degree $m$.
The \linebreak
$\glN$-module $L_{m}$ is irreducible with highest weight $(m,0,\dots,0)$ and
highest weight vector  $x_1^{m}$.

We consider  $\glN$ and $\glM$
simultaneously. To distinguish generators, modules, etc., we
indicated the dependence on $N$ and $M$ explicitly, for example, $E_{ij}\NNN\>$,
$L_m\MMM $.

Let
$P_{MN} = \C [x_{11},\dots,x_{M1},\dots,x_{1N},\dots,x_{MN}]$
be the space of polynomials of $MN$ variables.
We define the $\glM$-action on $P_{MN}$ by
$\displaystyle E_{ab}\MMM\,
\mapsto \sum_{i=1}^N\,x_{ai} \frac \partial{\partial x_{bi}}$
and  the $\glN$-action by
$\displaystyle E_{ij}\NNN\,
\mapsto \sum_{a=1}^M\,x_{ai} \frac \partial{\partial x_{aj}}.
$
There are two isomorphisms of vector spaces, 
\bean\label{Miso}
{}&&
\\
\bigl(\C[x_1\dots, x_M]\big)^{\ox N}\!\!\!\to P_{MN},
&& 
(p_1\ox\dots\ox p_N) (x_{11}, \dots,x_{MN})\mapsto 
\prod_{i=1}^N\,p_i(x_{1i},\dots,  x_{Mi}),
\notag
\\
\bigl(\C[x_1,\dots, x_N]\bigr)^{\ox M}\!\!\! \to P_{MN} ,
&&
(p_1\ox\dots\ox p_M)(x_{11},\dots,x_{MN})  \mapsto
\prod_{a=1}^M\,p_a(x_{a1},\dots, x_{aN}) .
\notag
\eean
Under these isomorphisms, $P_{MN}$ is isomorphic to $(\Vb\MMM) ^{\ox N}$ as a $\glM$-module
and  to $({\Vb}^{\langle N \rangle})^{\ox M}$ as a $\glN$-module.

Fix  $\bs n=(n_1,\dots,n_N)\in\Zp^{N}$ and
$\bs m=(m_1,\dots, m_M)\in\Zp^{M}$ with  $\sum_{i=1}^N n_i = \sum_{a=1}^M m_a$.
Isomorphisms \Ref{Miso}  induce an isomorphism of the weight subspaces,
\bean\label{duality isom}
 \LNM\ \simeq \LMN\ .
\eean
Under this isomorphism the KZ and dynamical 
Hamiltonians interchange,
\bea
H\NNN_{a}(\bs\la,\bs z)\,= \,G\MMM_{a}(\bs z,\bs \la)\,, 
 \qquad 
G\NNN_{i}(\bs\la,\bs z)\,=\,H\MMM_{i}(\bs z,\bs\la)\,, 
\eea
for $a=1,\dots, M$,\  $i=1,\dots,N$,\ \cite{TV}.

\bigskip

 Let $t^{\langle \bs n\rangle}$ be an admissible critical point of the master function
$\Phi(\,\cdot\,;\bs \la; \bs z; \bs m)$. Let 
$(V, \bs z)$ be the associated space of the  $(\bs \la, \bs z, \bs n, \bs m)$-type.
Let $(U, \bs \la)$ be its bispectral dual space of the  $(\bs z, \bs \la, \bs m, \bs n)$-type.
Assume that $(U, \bs \la)$ is admissible. Let 
$t^{\langle \bs m\rangle}$ be the associated admissible critical point of the master function
$\Phi(\,\cdot\,;\bs z; \bs \la; \bs n)$.

\bigskip

\noindent
\begin{conj}\label{conjecture}
 The corresponding Bethe vectors
$\omega(t^{\langle \bs n\rangle})\in \LMN$ and
$\omega(t^{\langle \bs m\rangle})\in \LNM$ are proportional under the duality isomorphism
\Ref{duality isom}.

\end{conj}

\section{The case $N =  M  = 2$}\label{sec m=n}

\subsection{Proof of Conjecture \ref{conjecture} for $N =  M  = 2$}\label{proof of conj}

Let $\bs n=(n_1,n_2)$ and $\bs m=(m_1,m_2)$ be two vectors of nonnegative integers such that
$n_1+n_2=m_1+m_2$. 

For $m\in\Zp$, let $L_m$ be the irreducible $\glt$-module with highest weight
$(m,0)$ and highest weight vector $v_m$. The vectors $E_{21}^iv_m, 
\ i=0,\dots, m$, form a basis in $L_m$.

Set 
$$
\alpha = \max\,(0, n_2-m_1)\ , \qquad
\beta =\min \,(m_2,n_2)\ .
$$
The vectors
\bea
\frac{E_{21}^{n_2-i}v_{m_1}}{(n_2-i)!}\ox \frac{E_{21}^iv_{m_2}}{i!}\ ,
\qquad
\alpha\leq i\leq \beta\ ,
\eea
form a basis in 
$\LMN\ =\ (L_{m_1}\ox L_{m_2})[n_1,n_2]$.
The vectors
\bea
 \frac{E_{21}^{m_2-i}v_{n_1}}{(m_2-i)!}\ox \frac{E_{21}^iv_{n_2}}{i!}\ ,
\qquad
\alpha  \leq i \leq \beta\ ,
\eea
form a basis in 
$\LNM\ =\ (L_{n_1}\ox L_{n_2})[m_1,m_2]$.
Isomorphism \Ref{duality isom} identifies the vectors with the same index $i$.

Notice that the map
\bean\label{weyl formula}
\frac{E_{21}^{n_2-i}v_{m_1}}{(n_2-i)!}\ox \frac{E_{21}^iv_{m_2}}{i!}
\ \mapsto
\frac{E_{21}^{m_1-n_2+i}v_{m_1}}{(m_1-n_2+i)!}\ox 
\frac{E_{21}^{m_2-i}v_{m_2}}{(m_2-i)!}
\eean
defines the Weyl isomorphism 
\bean\label{weyl isom}
(L_{m_1}\otimes L_{m_2})[n_1,n_2]\  \to (L_{m_1}\otimes L_{m_2})[n_2,n_1]\ .
\eean

\bigskip

Fix $\bs \la = (\la_1,\la_2)$ and $\bs z=(z_1,z_2)$ each 
 with distinct coordinates. 
{\it The universal rational 
$\LMN$-valued function} is the function
\bea
\omega (t_1, \dots , t_{n_2})\ =\!\!\! 
\sum_i
{C_i}\
 \frac{E_{21}^{n_2-i}v_{m_1}}{(n_2-i)!}\ \ox \frac{E_{21}^iv_{m_2}}{ i!}\ \ ,
\eea
where 
\bean\label{univ function}
C_i(t_1,\dots,t_{n_2})\ = \ \Sym\
 \prod_{j=1}^{n_2-i} 
\frac 1 {t_j - z_1}\
\prod_{j=1}^{i} 
\frac 1 {t_{n_2+j-i} - z_2}
\eean
and $\Sym \, f(t_1,\dots , t_{n})\ =\ \sum_{\sigma\in\Sigma_{n}} 
f(t_{\sigma(1)},\dots , t_{\sigma(n)})$, see \cite{FMTV}.

The universal rational function is
symmetric with respect to the group $\Sigma_{n_2}$ of permutations of 
variables $t_1,\dots,t_{n_2}$.
The $\Sigma_{n_2}$-orbit of a 
point $t^{\langle \bs 2 \rangle}=(t_1,\dots,t_{n_2})$ is represented by
the polynomial $p_2(x) = (x-t_1)\dots (x-t_{n_2})$.

\begin{lem}\label{lem y and c_i} Let
 $p_2(x) = (x-t_1)\dots (x-t_{n_2})$ be a polynomial. Then there
exist numbers $c_1, \dots , c_{n_2}$ 
such that
\bean\label{y and c_i}
p_2(x)\ =\ \sum_{i=1}^{n_2}\
c_i\
\frac{(x-z_1)^{n_2-i}}{(n_2-i)!}\ \frac{(x-z_2)^i}{i!}\ .
\eean
Moreover, 
 $c_i \ = \ (-1)^{i}\,(z_1-z_2)^{n_2}\ C_i$.
\hfill
$\square$
\end{lem}

Let $V = \langle p_1(x)e^{\lambda_1 x}, p_2(x)e^{\lambda_2 x}\rangle$.
Let $(V, \bs z)$ be a space of the $(\bs \la, \bs z, \bs n, \bs m)$-type. The
special fundamental operator of $(V, \bs z)$ has the form 
\bea
D &=& (x-z_1)(x-z_2) (\p - \lambda_1)(\p - \lambda_2) +
\phi _{11} (x-z_1) (\p - \lambda_1) +
\\
&&
\phi _{12} (x-z_1) (\p - \lambda_2) +
\phi _{21} (x-z_2) (\p - \lambda_1) +
\phi _{22} (x-z_2) (\p - \lambda_2) ,
\eea
where 
$\phi_{ij}$ are suitable numbers such that
\bea
\phi_{21}+\phi_{22}= - m_1  ,
& \qquad &
\phi_{11}+\phi_{12}= - m_2 ,
\\
\phi_{12}+\phi_{22}= - n_1 ,
& \qquad &
\phi_{11}+\phi_{21} = - n_2 .
\eea
Equation $D f = 0$ has a solution
\bean\label{sol 2}
f\ =\ e^{\lambda_2 x}\, p_2(x)\ =\ e^{\lambda_2 x} 
 \sum_{i}\ {}
c_i\
\frac{(x-z_1)^{n_2-i}}{(n_2-i)!}\ \frac{(x-z_2)^i}{i!}\ .
\eean
Substituting this expression to the differential equation we obtain  relations for
the coefficients,
\bean\label{rec eqns}
&&
(n_2-i)(m_2-i)\, c_{i+1}\ +\ i(n_1-m_2+i)\, c_{i-1}\ +
\\
&&
\phantom{aaa}
(\,-2i^2 + i(2n_2-m_1+m_2)-n_2m_2 +
(\lambda_1-\lambda_2)(z_1-z_2)(i+\phi_{11}))\ c_i\ =\ 0\ .
\notag
\eean
Notice that values $i=0, n_2-m_1, m_2, n_2$ are the values of $i$ for which
 equation \Ref{rec eqns} does not contain $c_{i-1}$ or $c_{i+1}$.

\begin{lem} 
Equations
\Ref{rec eqns} with $i$ such that
$\alpha \leq i\leq \beta $ form a closed
system of equations with respect to 
$c_j$ such that
$\alpha \leq j \leq \beta$. 
\hfill
$\square$
\end{lem}

Equations \Ref{rec eqns} have two symmetries.

The first symmetry has the following form. Replace in
\Ref{rec eqns} the parameters $\la_1,$ $ \la_2,$ $ z_1,$ 
$ z_2,$ $ n_1, $ $n_2, m_1, m_2$
by $ z_1, z_2, \la_1, \la_2,  m_1, m_2, n_1, n_2$, respectively.
Then \Ref{rec eqns} does not change.

The second symmetry has the following form. In the coefficients of
\Ref{rec eqns} replace the parameters $\la_1, \la_2,  n_1, n_2,  \phi_{11}, i$
by $ \la_2, \la_1,  n_2, n_1,  \phi_{12}, m_2-j$, respectively.
Then \Ref{rec eqns} takes the form
\bean\label{rec eqns dual}
&&
(n_2-j)(m_2-j)\, c_{i-1}\ +\ j(n_1-m_2+j)\, c_{i+1}\ +
\\
&&
\phantom{aaa}
(\,-2j^2 + j(2n_2-m_1+m_2)-n_2m_2 +
(\lambda_1-\lambda_2)(z_1-z_2)(j+\phi_{11}))\ c_i\ =\ 0\ .
\notag
\eean

The symmetries imply
\begin{theorem}\label{bethe thm}
Let $V = \langle p_1(x)e^{\lambda_1 x}, p_2(x)e^{\lambda_2 x} \rangle$
and
$U=\langle q_1(u)e^{z_1 u}, q_2(u)e^{z_2 u} \rangle$.
Let $(V, \bs z)$ be a space of the $(\bs \la, \bs z, \bs n, \bs m)$-type and 
$(U, \bs \la)$   the special bispectral dual space
of the $(\bs z, \bs \la, \bs m, \bs n)$-type. Write
\bea
 e^{\lambda_2 x}\, p_2(x)\ =\ e^{\lambda_2 x} 
 \sum_{i=0}^{n_2}\ {}
c_i\
\frac{(x-z_1)^{n_2-i}}{(n_2-i)!}\ \frac{(x-z_2)^i}{i!}\ ,
\\
e^{\lambda_1 x}\, p_1(x)\ =\ e^{\lambda_1 x} 
 \sum_{i=0}^{n_1}\ {}
d_i\
\frac{(x-z_1)^{n_1-i}}{(n_1-i)!}\ \frac{(x-z_2)^i}{i!}\ ,
\\
e^{z_2 u}\, q_2(u)\ =\ e^{z_2 u} 
 \sum_{i=0}^{m_2}\ {}
e_i\
\frac{(u-\la_1)^{m_2-i}}{(m_2-i)!}\ \frac{(u-\la_2)^i}{i!}\ .
\eea
Then $c_i=d_{m_2-i}=e_i$ for $\alpha\leq i\leq \beta$.
\end{theorem}
\begin{cor}\label{cor dual}
The Bethe vectors 
\bea
&&
\sum_{i=\alpha}^{\beta}
{c_i}\
 \frac{E_{21}^{n_2-i}v_{m_1}}{(n_2-i)!}\ \ox \frac{E_{21}^iv_{m_2}}{ i!}\  ,
\qquad 
\sum_{i=\alpha}^{\beta}
{d_i}\
 \frac{E_{21}^{n_1-i}v_{m_1}}{(n_1-i)!}\ \ox \frac{E_{21}^iv_{m_2}}{ i!}\  ,
\\
&&
\phantom{aaaaaaaaa}
\sum_{i=\alpha}^{\beta}
{e_i}\
 \frac{E_{21}^{m_2-i}v_{n_1}}{(m_2-i)!}\ \ox \frac{E_{21}^iv_{n_2}}{ i!} 
\eea
are identified by isomorphisms  \Ref{weyl formula} and \Ref{duality isom}.
\end{cor}
The fact that the first and second Bethe vectors are identified by the Weyl isomorphism
is proved also in \cite{MV3} in a different way.

\subsection{Critical points of master functions for $N =  M  = 2$}\label{sec crit points}
Under conditions of Section \ref{proof of conj}
consider the three associated master functions,
\bea
&&
\Phi (t_1,\dots,t_{n_2}; \bs\la; \bs z; \bs m) =
\exp\left(\la_1(m_1z_1 +m_2z_2) + (\la_2-\la_1)\sum_{j=1}^{n_2} t_j \right) \times
\\
&&
\phantom{aaaaaaaaaaa}
(z_1-z_2)^{m_1m_2}\ 
\prod_{j=1}^{n_{2}}
(t_j-z_1)^{-m_1}(t_j-z_2)^{-m_2}\!\!\!
\prod_{1\leq j<j'\leq n_2}\!\!\!
(t_j-t_{j'})^{2}\ ,
\\
&&
\Phi (t_1,\dots,t_{n_1}; \la_2,\la_1; \bs z; \bs m) =
\exp\left(\la_2(m_1z_1 +m_2z_2) + (\la_1-\la_2)\sum_{j=1}^{n_1} t_j \right) \times
\\
&&
\phantom{aaaaaaaaaaa}
(z_1-z_2)^{m_1m_2}\ 
\prod_{j=1}^{n_{1}}
(t_j-z_1)^{-m_1}(t_j-z_2)^{-m_2}\!\!\!
\prod_{1\leq j<j'\leq n_1}\!\!\!
(t_j-t_{j'})^{2}\ .
\\
&&
\Phi (t_1,\dots,t_{m_2}; \bs z; \bs \la; \bs n) =
\exp\left(z_1(n_1\lambda_1 +n_2\lambda_2) + (z_2-z_1)\sum_{j=1}^{m_2} t_j \right) \times
\\
&&
\phantom{aaaaaaaaaaa}
(\la_1-\la_2)^{n_1n_2}\ 
\prod_{j=1}^{m_{2}}
(t_j-\la_1)^{-n_1}(t_j-\la_2)^{-n_2}\!\!\!
\prod_{1\leq j<j'\leq m_2}\!\!\!
(t_j-t_{j'})^{2}\ .
\eea
The functions are symmetric with respect to permutations of coordinates.

Set\
$d \ =\  \dim\ (L_{m_1}\otimes L_{m_2})[n_1,n_2] =
\dim \ (L_{m_1}\otimes L_{m_2})[n_2,n_1] = 
\dim\ (L_{n_1}\otimes L_{n_2})[m_1,m_2].$ 
By \cite{MV3} for generic $\bs \la$ and $\bs z$
each of the three functions has exactly $d$ orbits of critical points 
and each of the orbits consists of non-degenerate critical points.

\begin{theorem} Let $\bs \la$ and $\bs z$ be generic.
Let $(t_1, \dots , t_{n_2})$ be an admissible critical point of 
$\Phi (\,\cdot\,; \bs\la; \bs z; \bs m)$. Let
$(V, \bs z)$ be the space of the $(\bs \la, \bs z, \bs n, \bs m)$-type associated with the
orbit of the critical point. 
 Let $(U, \bs \la)$  be the special bispectral dual space to $V$.
Let $U=\langle q_1(u)e^{z_1 u}, q_2(u)e^{z_2 u} \rangle$.
Then the polynomial
$p_1$ represents the orbit of an admissible critical point
of $\Phi (t_1,\dots,t_{n_1}; \la_2,\la_1; \bs z; \bs m)$
and the polynomial $q_2$ represents the orbit of an admissible critical point
of $\Phi (t_1,\dots,t_{m_2}; \bs z; \bs \la; \bs n)$.
\hfill
$\square$
\end{theorem}

\bigskip

Let
$F(t_1,\dots , t_k; \la_1,\la_2,z_1,z_2)$ be a function
of variables $t_1,\dots,t_k$ depending on parameters
$\la_1,\la_2,z_1,z_2$. Let $C \subset \C^k\times\C^4$ be the subset of critical points
of $F$ with respect to variables $t_1,\dots,t_k$. Consider the space $\C^8$ with coordinates
$\la_1,\la_2,z_1,z_2, \la_1^*,\la_2^*,z_1^*,z_2^*$ and symplectic form
$d\la_1\wedge d\la_1^* + \la_2\wedge
d\la_2^*+dz_1\wedge
dz_1^* +dz_2\wedge dz_2^*$.
Define a map $C \to \C^8$ by the formula
\bea
(t; \bs \la;\bs z )\
\mapsto\
(\bs \la; \bs z;
\frac{\partial F(t; \bs \la;\bs z )}{\partial \la_1},
\frac{\partial F(t; \bs \la;\bs z )}{\partial \la_2},
\frac{\partial F(t; \bs \la;\bs z )}{\partial z_1},
\frac{\partial F(t; \bs \la;\bs z )}{\partial z_2} ) \ .
\eea
The image is a Lagrange variety.

\begin{theorem}
The Lagrange varieties corresponding to the logarithms of the three master functions coincide.
\end{theorem}

\begin{proof}
The three Bethe vectors in Corollary \ref{cor dual} are identified by
isomorphisms \Ref{weyl formula} and \Ref{duality isom}. The KZ and dynamical Hamiltonians
acting in $(L_{m_1}\otimes L_{m_2})[n_1,n_2]$, 
$(L_{m_1}\otimes L_{m_2})[n_2,n_1]$, and $(L_{n_1}\otimes L_{n_2})[m_1,m_2]$  
also are identified by 
isomorphisms \Ref{weyl formula} and \Ref{duality isom}. Hence the three Bethe vectors
have the same eigenvalues with respect to the Hamiltonians. But the eigenvalues are given
by partial derivatives of the logarithms of the corresponding master functions, see Theorem
\ref{thm RV}.  This proves that the Lagrange varieties coincide.
\end{proof}

\section{Baker-Akhieser functions and bispectral involution}\label{baker}
For $\lambda \in \C$, we call a vector subspace $W_\la \subset \C[t]$ 
 {\it admissible at $\la$} if there exists $m\in\Z_{\geq 0}$ such that
$(t-\la)^m\,\C[x] \subset W_\la$ and there exists $f\in W_\lambda$ such that $f(\lambda)\neq 0$. 
We call a vector
subspace $W \subset \C[t]$ {\it  admissible} if 
\bea
W\ =\ \cap_{i=1}^n\ W_{\la_i}\ 
\eea
where $n$ is a natural number, $\la_1,\dots,\la_n\in\C$ are distinct complex numbers,
and for any $i$, $W_{\la_i}$ is admissible at $\la_i$.

We denote the set of all admissible subspaces by $\Gr$ and call {\it the Grassmannian of
admissible subspaces}. The Grassmannian parameterizes the rational solutions 
(vanishing at infinity) to the KP hierarchy, see \cite{W}.

Each $W_{\lambda_i}$ can be defined by a finite set of linear equations on the Taylor coefficients
at $\la_i$. Namely, there exist a positive integer $N_i$, a sequence of
 integers $0 < n_{i1}<\dots < n_{iN_i}$ and complex numbers
$c_{i,j,a}$, $a=0,\dots, n_{ij}$, such that $c_{i,j,n_{ij}} \neq 0$ for all $i,j$ and
\bean\label{adm subspace}
W_{\la_i}\ =\ \{\ r \in \C[t]\ \ |\ 
\ \sum_{a=1}^{n_{ij}}\ c_{i,j,a}\ r^{(a)}(\lambda_i) = 0 ,\
 j = 1, \dots , N_i\ \}\ .
\eean

For an admissible subspace $W$ define the complex vector space $V$ as
the space spanned by functions $p_{ij} e^{\la_i x}$, $i=1,\dots , n$,\
$j=1,\dots , N_i$, where $p_{ij} = \sum_{a=1}^{n_{ij}}\ c_{ija} x^a
e^{\la_i x}$.  The space $V$ is a space of quasi-polynomials of type
considered in Section \ref{subs exponentials}.

On the other hand, having a space $V$ of type considered in Section
\ref{subs exponentials} one can recover an admissible subspace
$W\subset \C[t]$ by formula \Ref{adm subspace}.

\bigskip

Let $W\subset\C[t]$ be an admissible subspace and $V$ the associated
space of quasi-polynomials.  Define the algebra \ $A = \{ p\in \C[t]\
| \ p(t)\,W \subset W \}.$ An equivalent definition is $A = \{ p\in
\C[t]\ | \ p(\partial_x)\,V \subset V \}.$

Let $\bar D_V$ be the monic fundamental differential operator
of $V$. 
The function
\bea
\Psi_W (x,\xi) \ = \ {\prod_{i=1}^n(\xi-\la_i)^{-N_i}} \
\bar D_V\ e^{x\xi} \ 
\eea
is called {\it
the stationary Baker-Akhieser function 
of the admissible space $W$}. Introduce the rational
function $\psi_W (x,\xi)$ by the formula $\Psi_W (x,\xi) = \psi_W (x,\xi) e^{x\xi}$.
The function $\psi_W(x,\xi)$ expands as a power series in $x^{-1}, \xi^{-1}$
of the form
$$
\psi_W(x,\xi) = 1 + \sum_{i,j=1}^\infty c_{ij} x^{-i} \xi^{-j}.
$$
It is easy to see that for every $p\in A$, there exists a linear
differential operator $L_p(x,\p)$ such that $L_p(x,\p) \bar D_V = \bar
D_V p(\p)$. As a corollary one concludes that 
$$
L_p(x,\p) \,\Psi_W(x,\xi)\ =\ p(\xi)\, \Psi_W(x,\xi) \ .
$$

There is the bispectral involution for points of the Grassmanian of
admissible subspaces.  It is described in terms of the Baker-Akhiezer
functions. An admissible space $W$ is mapped to an admissible space $Y$ if
$\Psi_Y(x,\xi)=\Psi_W(\xi,x)$, see \cite{W}.

\begin{theorem}
Let $V$ and $U$ be 
two spaces of quasi-polynomials bispectral dual with respect to the integral transform 
of Section \ref{sec int transform}. 
Then the corresponding admissible spaces \linebreak
$W, Y \in \Gr$ are bispectral dual
with respect to their Baker-Akhieser functions.
\end{theorem}

The theorem follows from part (iii) of Theorem \ref{thm int transform}.

The points $W$ of the  Grassmanian $\Gr$ 
are in bijection with classes of pairs of matrices $Z, \Lambda$ such that
$[Z,\Lambda]+I$ is a rank one matrix \cite{W}. Then
$$
\psi_W(x,\xi)\ = \ {\rm det}( I + (x-Z)^{-1} (\xi-\Lambda)^{-1} )\ .
$$
The bispectral involution corresponds to the involution $(Z,\Lambda) \mapsto (\Lambda^*,Z^*)$.

The simplest example of both diagonalizable $Z$ and $\Lambda$ corresponds
to the case of the $(\glN, \glN)$ duality and the weight subspace
$({\C}^N)^{\otimes N}[1,...,1]$
that is, the weight subspace of $\slN$ weight zero. This is a remarkable weight
subspace; it is naturally isomorphic to the group algebra of the symmetric
group $\Sigma_N$.

\section{Appendix: Proof of Proposition \ref{lem on multiplicity}}\label{appendix}

Let $\g$ be a simple Lie algebra. Fix its Gauss decomposition
$\g=\n^+\oplus\h\oplus\n^-$. For a $\g$-module $U$ let $U[\lambda]$ be
the weight subspace of weight $\lambda$ and let $U^{\n^+}[\lambda]$ be
the subspace of singular vectors of weight $\lambda$. Let $M_\lambda$ be
the Verma module over $\g$ of highest weight $\lambda$, and let $L_\lambda$ be
the irreducible $\g$-module of highest weight $\lambda$.
\begin{lem}
\label{HomM}
Let $U$ be a $\g$-module. Then
$\dim\bigl(\Hom_\g(M_\lambda,U)\bigr)=\dim\bigl(U^{\n^+}[\lambda]\bigr)$.
\hfill
$\square$
\end{lem}
\begin{lem}\label{HomL}
Let $U$ be a finite-dimensional $\g$-module. Then
$$
\dim\bigl(\Hom_\g(L_\lambda,U)\bigr)=\dim\bigl(U^{\n^+}[\lambda]\bigr)\,.
$$
\end{lem}
\begin{proof}
Since $U$ is finite-dimensional and the only possible finite-dimensional
quotient of $M_\lambda$ is $L_\lambda$, any homomorphism
${\phi\in\Hom_\g(M_\lambda,U)}$ factors through $L_\lambda$.
Therefore, the canonical embedding
$\Hom_\g(L_\lambda,U)\hookrightarrow\Hom_\g(M_\lambda,U)$ is an isomorphism,
and the claim follows from Lemma \ref{HomM}.
\end{proof}
For a $\g$-module $U$, the dual $\g$-module $U^*$ is defined by the rule
$x|_{U^*}=(-x|_U)^*$ for any $x\in\g$. Let $w$ be the longest element
of the Weyl group of $\g$.
\begin{lem}
Let $\lambda$ be a dominant integral weight. Then the $\g$-module $L_\lambda^*$
dual to $L_\lambda$ is isomorphic to $L_{-w(\lambda)}$.
\hfill
$\square$
\end{lem}
\begin{lem}
\label{Using}
Let $\lambda$ be a dominant integral weight, and let $U$ be
a finite-dimensional $\g$-module. The multiplicity of the trivial $\g$-module
in the tensor product $U\ox L_\lambda$ equals
$\dim\bigl(U^{\n^+}[-w(\lambda)]\bigr)$.
\end{lem}
\begin{proof}
The multiplicity of the trivial $\g$-module in the tensor product
$U\ox L_\lambda$ equals $\dim\bigl(\Hom_\g(\C,U\ox L_\lambda)\bigr)$.
Since $\Hom_\g(\C,U\ox L_\lambda)=\Hom_\g(L_\lambda^*,U)=
\Hom_\g(L_{-w(\lambda},U)$, the claim follows from Lemma \ref{HomL}.
\end{proof}
Let $r=\rank\,\g$. Let $(\,,\,)$ be the invariant bilinear form on $\h^*$,
and let $\alpha_1^\vee,\ldots\alpha_r^\vee$ be simple coroots.
Let $e_1,\ldots e_r\in\n^+$ and $f_1,\ldots f_r\in\n^-$ be the Chevalley
generators corresponding to the simple roots.
\begin{lem}
\cite[7.2.5]{Dix}
\label{Dix725}
Let $\lambda$ be a dominant integral weight. Then
$$
L_\lambda=M_\lambda/\{f_1^{(\lambda,\alpha_1^\vee)+1}v_\lambda,\ldots,
f_r^{(\lambda,\alpha_r^\vee)+1}v_\lambda\}\,,
$$
where $v_\lambda$ is the highest weight vector in $M_\lambda$.
\hfill
$\square$
\end{lem}
Let $\varpi:\g\to\g$ be the involution such that $\varpi(e_i)=-f_i$,
\ $i=1,\ldots,r$. For a $\g$-module $U$ define the $\g$-module $\overline U$
on the same vector space by the rule $x|_{\overline U}=\varpi(x)|_U$ for any
$x\in\g$.
\begin{lem}
Let $\lambda$ be a dominant integral weight. The $\g$-module
$\overline{L_\lambda^*}$ is isomorphic to $L_\lambda$.
\hfill
$\square$
\end{lem}

Let $\mathfrak{b}^\pm$ 
$=\h\oplus\n^\pm$. For a weight $\mu$ let $\C_\mu$ be
the one-dimensional $\h$-module such that $x\in\h$ acts as $\mu(x)$. We treat
$\C_\mu$ also as a module over $\mathfrak{b}^+$ or $\mathfrak{b}^-$ in which
the subalgebras $\n^\pm$
act trivially.
\begin{lem}
\label{HomMb}
Let $V$ be a $\g$-module. Let $v_\lambda$ be a highest weight vector of
the Verma module $M_\lambda$ over $\g$. Then the map
$\Hom_{\mathfrak{b}^-}(\C_\mu\ox M_\lambda,V)\to V[\lambda+\mu]$,
\ $\phi\mapsto\phi(1\ox v_\lambda)$ is a bijection.
\hfill
$\square$
\end{lem}

\begin{lem}\label{Ue}
Let $\lambda\,,\,\mu$ be dominant integral weights.
Let $U$ be a finite-dimensional \linebreak
$\g$-module. Then
$$
\dim\bigl((U\ox L_\lambda)^{\n^+}[\mu]\bigr)=
\dim\bigl\{u\in U[\mu-\lambda]\ \big|\ e_1^{(\lambda,\alpha_1^\vee)+1}u=0,
\ldots,e_r^{(\lambda,\alpha_r^\vee)+1}u=0\bigr\}\,.
$$
\end{lem}
\begin{proof}
We have
$$
(U\ox L_\lambda)^{\n^+}[\mu]=
\Hom_{\mathfrak{b}^+}\bigl(\C_\mu,U\ox L_\lambda\bigr)=
\Hom_{\mathfrak{b}^+}\bigl(\C_\mu\ox L_\lambda^*,U\bigr)={}
$$
$$
{}=\Hom_{\mathfrak{b}^-}\bigl(\overline{\C_\mu}\ox\overline{L_\lambda^*},
\overline U\bigr)=
\Hom_{\mathfrak{b}^-}\bigl(\C_{-\mu}\ox L_\lambda,\overline U\bigr)\,.
$$
Let $v_\lambda$ be a highest weight vector of $L_\lambda$.
By Lemma \ref{Dix725} and Lemma \ref{HomMb} the map
$$
\Hom_{\mathfrak{b}^-}(\C_{-\mu}\ox L_\lambda,\overline U)\to
\bigl\{\bar u\in\overline U[\lambda-\mu]\ |
\ f_1^{(\lambda,\alpha_1^\vee)+1}\bar u=0,
\ldots,f_r^{(\lambda,\alpha_r^\vee)+1}\bar u=0\bigr\}\,,
$$
$$
\phi\mapsto\phi(1\ox v_\lambda)
$$
is a bijection. At last, observe that
$$
\bigl\{\bar u\in\overline U[\lambda-\mu]\ |
\ f_1^{(\lambda,\alpha_1^\vee)+1}\bar u=0,
\ldots,f_r^{(\lambda,\alpha_r^\vee)+1}\bar u=0\bigr\}={}
$$
$$
{}=\bigl\{u\in U[\mu-\lambda]\ \big|\ e_1^{(\lambda,\alpha_1^\vee)+1}u=0,
\ldots,e_r^{(\lambda,\alpha_r^\vee)+1}u=0\bigr\}\,.
$$
To avoid confusion, notice that though $U[\mu-\lambda]$ and
$\overline U[\lambda-\mu]$ coincide as vector spaces, they are considered
as subspaces of $\g$-modules $U$ and $\overline U$, respectively, and
the corresponding $\g$-actions do not coincide.
\end{proof}
\begin{prop}
\label{gtriv}
Let $\lambda\,,\,\mu$ be dominant integral weights. Let $U$ be 
a finite-dimensional $\g$-module. Then the multiplicity of the trivial
$\g$-module in the tensor product $U\ox L_\lambda\ox L_\mu$ equals
$$
\dim\bigl\{u\in U[-\lambda-w(\mu)]\ \big|\ e_1^{(\lambda,\alpha_1^\vee)+1}u=0,
\ldots,e_r^{(\lambda,\alpha_r^\vee)+1}u=0\bigr\}\,.
$$
\end{prop}
\noindent
The statement follows from Lemmas \ref{Using} and \ref{Ue}.

\bigskip
Return to the situation of Section \ref{number}.  
For
a ramification sequence $\bs c$ denote by $L_{\bs c}$ the finite-dimensional
irreducible $\glN$-module with highest weight $\bs c$. Any $\glN$-module $U$
has a natural structure of an $\slN$-module denoted by $\widetilde U$.

A $\glN$-module $U$ is called polynomial of degree $d$ if it is a submodule
of the tensor product $(\C^N)^{\ox d}$ of the vector representations of $\glN$.
For any ramification sequence $\bs c$ the module $L_{\bs c}$ is polynomial.
\begin{lem}
\label{Uweights}
Let $U$ be a polynomial $\glN$-module of degree $d$. If a weight subspace
$U[n_1,\ldots,n_N]$ is nonzero, then $n_1,\ldots n_N$ are nonnegative integers
and $\sum_i n_i=d$.
\hfill
$\square$
\end{lem}

Let $U$ be a polynomial $\glN$-module of degree $d$ and let $n_1,\dots, n_N$ be
nonegative integers such that $\sum_i n_i=d$. Define the ramification sequences
\bea
\bs c(\infty_1)& =& (n_2+\dots+n_N, n_3+\dots+ n_N,\dots,n_{N},0)\ ,
\\
\bs c(\infty_2) & = & (n_1+\dots+n_{N-1}, n_1+\dots+ n_{N-2},\dots,n_{1},0)\ .
\eea

\begin{prop}\label{lemmma on multiplicity}
The multiplicity of the trivial $\slN$-module in the tensor product of
$\slN$-modules $\widetilde U\ox \widetilde L_{\bs c(\infty_1)}\ox
\widetilde L_{\bs c(\infty_2)}$ equals $\dim U[n_1,\dots,n_N]$.
\end{prop}
\begin{proof}
Let $\lambda\,,\,\mu$ be the $\slN$ weights corresponding to the $\glN$ weights
$\bs c(\infty_1)\,,\bs c(\infty_2)$, respectively.
Then $\widetilde U[-\lambda-w(\mu)]=U[n_1,\dots,n_N]$ and
$(\lambda,\alpha_i^\vee)=n_{i+1}$, \ $i=1,\ldots,N-1$.
By Lemma \ref{Uweights} we have  $e^{n_{i+1}+1}u=0$, \ $i=1\ldots,N-1$,
for any $u\in U[n_1,\dots,n_N]$. Hence, the required statement follows from
Proposition \ref{gtriv}.
\end{proof}

\bigskip

\end{document}